\documentclass[11pt,oneside,leqno]{article}
\usepackage[utf8]{inputenc}
\usepackage[russian]{babel}
\usepackage{amssymb,amsmath,amsthm}
\renewcommand{\subsection}{\refstepcounter{subsection}%
\par\bigskip\noindent\textbf{\upshape\thesubsection. }}
\renewcommand{\subsubsection}{\refstepcounter{subsubsection}%
\par\medskip\noindent\textbf{\upshape\thesubsubsection.  }}
\renewcommand{\paragraph}{\refstepcounter{paragraph}%
\par\smallskip\noindent\textbf{\upshape\theparagraph. }}

\numberwithin{equation}{subsection}

\renewcommand{\thesubsection}{\arabic{subsection}}

\renewcommand{\land}{\mathbin{\&}}
\newcommand{\kom}{\mathfrak c}

\title{\begin{flushleft}\normalsize УДК~510.25\end{flushleft}
Об одной теореме непрерывности для конструктивных функций}
\author{А.~А.~Владимиров\footnote{Работа поддержана РФФИ, грант \No~06-06-80326.}}
\date{}
\begin{document}
\renewcommand{\proofname}{{\upshape Д\,о\,к\,а\,з\,а\,т\,е\,л\,ь\,с\,т\,в\,о.}}
\maketitle
\begin{abstract}
В статье устанавливается факт локальной равномерной непрерывности всюду
определённого и сохраняющего предкомпактность подмножеств конструктивного
отображения полного метрического пространства в полное метрическое пространство.
Этот результат может быть рассмотрен в качестве конструктивного уточнения
интуиционистской теоремы Л.~Э.~Я.~Брауэра о веере.
\end{abstract}

\section{Введение}\label{par:1}
\subsection
Одним из центральных положений интуиционистской математики является теорема
Л.~Э.~Я.~Брауэра о веере (см., например, \cite[3.4.2]{Heyt:1965}), утверждающая
локальную равномерную непрерывность произвольного отображения полного
метрического пространства в метрическое пространство. Формулировка этой теоремы
опирается на не имеющее конструктивного истолкования понятие \emph{свободно
становящейся последовательности}, причём попытка его механической замены
понятием конструктивной последовательности приводит к утверждениям, легко
опровергаемым на контрпримерах (см., например \cite[комментарий~19]{Heyt:1965}).
В связи с этим высказывалось мнение, что для конструктивного направления в математике
теорема о веере может играть лишь эвристическую роль (см., например,
\cite[28]{ML:1975}). Аналогичной позиции придерживается школа Э.~Бишопа, согласно
представлениям которой локальная равномерная непрерывность рассматривается
как одна из составных частей определения понятия непрерывного отображения
(см., например, \cite[Ch.~4, Sect.~6]{BiBr:1985}).

Тем не менее, в рамках конструктивного направления в математике также оказывается
возможной характеризация класса локально равномерно непрерывных отображений
полных метрических пространств, не связанная с предварительными предположениями
о каких-либо континуальных свойствах таких отображений.  А именно, для локальной
равномерной непрерывности всюду определённого конструктивного отображения
\(f:X\to Y\) полного метрического пространства \(X\) в полное метрическое
пространство \(Y\) необходимо и достаточно, чтобы \(f\)-образом любого
предкомпактного подмножества также являлось некоторое предкомпактное
подмножество. В установлении этого факта и заключается основная цель
настоящей статьи.

\subsection
В дальнейшем мы будем придерживаться конструктивного понимания математических
суждений (см., например, \cite{Shan:1958}). В частности, утверждения вида
\begin{quote}
"<для любого конструктивного объекта \(x\) со свойством \(\mathcal P\) осуществим
конструктивный объект \(y\), находящийся в отношении \(\mathcal R\) к исходному
объекту \(x\)">
\end{quote}
будут рассматриваться нами как сокращённые варианты записи утверждений вида
\begin{quote}
"<осуществим алгорифм \(\mathfrak A\), применимый к любому конструктивному объекту
\(x\) со свойством \(\mathcal P\) и перерабатывающий его в конструктивный объект
\(y\), находящийся в отношении \(\mathcal R\) к исходному объекту \(x\)">.
\end{quote}
В качестве примера логико-математического языка, семантика которого построена
в соответствии с такими принципами, может быть указан язык \(\text{Я}_{\omega |}\)
ступенчатой семантической системы А.~А.~Маркова (см. \cite{Mark:1974}).

\subsection
Термин "<множество"> мы далее будем рассматривать как сокращение для термина
"<однопараметрическая формула логико-математического языка"> (см., например,
\cite[\S~7]{Shan:1958}, \cite[Введение, 7]{Ku:1973}). Знак равенства "<\(=\)">,
связывающий два множества, мы будем истолковывать как указание на равнообъёмность
этих множеств. Для записи отношения принадлежности конструктивного объекта множеству
и отношения вложения одного множества в другое мы будем использовать стандартные
теоретико-множественные символы "<\(\in\)"> и "<\(\subseteq\)"> (см., например,
\cite[Введение, 7]{Ku:1973}).

\subsection
Термины "<вещественное число"> и "<метрическое пространство"> мы далее будем
использовать в качестве сокращений для терминов "<конструктивное вещественное
число"> и "<конструктивное метрическое пространство с непустым носителем">.
Определения этих понятий могут быть найдены в \cite[Гл.~2, \S~2]{Ku:1973}
и \cite[Гл.~9, \S~1]{Ku:1973}, соответственно.

Для обозначения \emph{носителя} метрического пространства \(X\) мы будем
использовать запись \(\mathfrak S_X\), а для обозначения \emph{метрического
алгорифма} этого пространства "--- запись \(\rho_X\).

Термин "<предельная точка подмножества \(\mathfrak A\subseteq\mathfrak S_X\)">
мы будем связывать с точками \(x\in\mathfrak S_X\), внутри любой окрестности
которых квазиосуществима некоторая точка \(x'\in\mathfrak A\) (см., например,
\cite[Гл.~9, \S~3, Определение~4]{Ku:1973}). Соответственно, термин "<замкнутое
подмножество"> мы будем связывать с подмножествами, содержащими все свои предельные
точки.

Знак равенства "<\(=\)">, связывающий две точки некоторого метрического
пространства, мы будем истолковывать как указание на обращение в нуль расстояния
между этими точками.

\subsection
Через \(\mathbb R\) и \(\mathbb R^+\) мы будем обозначать, соответственно, множество
вещественных чисел и множество положительных вещественных чисел. Через \(\mathbf R\)
мы будем обозначать метрическое пространство с носителем \(\mathbb R\) и обычным
метрическим алгорифмом \(\rho_{\mathbf R}(x,x')\rightleftharpoons |x-x'|\).

\subsection
При ссылках на разделы статьи, не принадлежащие параграфу, внутри которого даётся
ссылка, мы будем дополнительно указывать номер параграфа. Аналогично, при ссылках
на формулы, не принадлежащие пункту, внутри которого даётся ссылка, мы будем
дополнительно указывать номер пункта.


\section{Пространство непустых компактов}\label{par:2}
\subsection\label{pt:21}
Пусть \(X\) "--- полное метрическое пространство. Для произвольных двух непустых
списков \(\zeta\rightleftharpoons\{\zeta_k\}_{k=1}^n\)
и \(\eta\rightleftharpoons\{\eta_l\}_{l=1}^m\) точек пространства \(X\)
может быть конструктивно определено \emph{расстояние Хаусдорфа}
\begin{equation}\label{eq:111}
	h_X(\zeta,\eta)\rightleftharpoons
	\sup\bigl(\sup\limits_{k\in [1,n]}
	\mathop{\inf\vphantom{p}}\limits_{l\in [1,m]}\rho_X(\zeta_k,\eta_l),
	\sup\limits_{l\in [1,m]}
	\mathop{\inf\vphantom{p}}\limits_{k\in [1,n]}\rho_X(\eta_l,\zeta_k)\bigr).
\end{equation}
Результат конструктивного пополнения множества всевозможных непустых списков точек
пространства \(X\) по расстоянию \eqref{eq:111} мы в дальнейшем будем называть
\emph{пространством непустых компактов} пространства \(X\) и обозначать через
\(\kom(X)\). Такой способ введения понятия компакта в рамках конструктивного
математического анализа указан, в частности, в работе \cite[4.5]{Shan:2000}.
Целью настоящего параграфа является описание некоторых свойств пространства
\(\kom(X)\), которые будут использованы нами при доказательстве основного результата
статьи. Приводимые в настоящем параграфе утверждения не претендуют в полной мере
на новизну (в частности, утверждение \ref{prop:135} в несколько других терминах
упоминается уже в работе \cite[\S~11.2]{Shan:1962}).

\subsection\label{pt:22}
Рассмотрим алгорифм \(c_X\), сопоставляющий произвольной паре непустых списков
\(\zeta\rightleftharpoons\{\zeta_k\}_{k=1}^n\) и \(\eta\rightleftharpoons
\{\eta_k\}_{k=1}^m\) точек пространства \(X\) список \(\theta\rightleftharpoons
\{\theta_k\}_{k=1}^{n+m}\) вида
\begin{equation}\label{eq:001}
	\theta_k\eqcirc\left\{\begin{aligned}
		&\zeta_k&\text{при }k\leqslant n,\\
		&\eta_{k-n}&\text{при }k>n,
	\end{aligned}\right.
\end{equation}
т.~е. результат конкатенации списков \(\zeta\) и \(\eta\). Для любых четырёх
непустых списков \(\zeta\), \(\eta\), \(\zeta'\) и \(\eta'\) точек пространства
\(X\) с очевидностью выполняется соотношение
\begin{flalign}
	\label{eq:003}
	&& h_X(c_X(\zeta,\eta),c_X(\zeta',\eta'))&\leqslant
	\sup\bigl(h_X(\zeta,\zeta'),\,h_X(\eta,\eta')\bigr).&
	\text{[\eqref{eq:001}, \ref{pt:21}\,\eqref{eq:111}]}&
\end{flalign}
Тем самым, алгорифм \(c_X\) допускает продолжение по непрерывности (см., например,
\cite[\S~10.2]{Shan:1962}) до конструктивной бинарной операции
\(\cup:\kom(X)\times \kom(X)\to \kom(X)\) \emph{объединения} непустых компактов. При этом
на основе предельного перехода легко устанавливается справедливость следующих
утверждений:

\subsubsection\label{prop:109}
{\itshape Операция объединения непустых компактов является сочетательной
и перестановочной} [\eqref{eq:001}, \ref{pt:21}\,\eqref{eq:111}].

\subsubsection\label{prop:110}
{\itshape Пусть \(K_1\), \(K_2\), \(K'_1\) и \(K'_2\) "--- четыре точки
пространства \(\kom(X)\). Тогда выполняется соотношение
\begin{flalign*}
	&& \rho_{\kom(X)}(K_1\cup K_2,K'_1\cup K'_2)&\leqslant
	\sup\bigl(\rho_{\kom(X)}(K_1,K'_1),\,\rho_{\kom(X)}(K_2,K'_2)\bigr).&
	\text{\upshape [\eqref{eq:003}]}&
\end{flalign*}
}

\subsection
Компакт \(K\in\mathfrak S_{\kom(X)}\) мы будем называть \emph{мажорирующим} компакт
\(K'\in\mathfrak S_{\kom(X)}\), если выполняется равенство \(K'\cup K=K\). Утверждение
о мажорировании компакта \(K'\) компактом \(K\) мы будем сокращённо записывать
в виде \(K'\leqslant K\).

Точку \(x\in\mathfrak S_X\) мы будем называть \emph{инцидентной} компакту
\(K\in\mathfrak S_{\kom(X)}\), если одночленный список \(\{x\}\) мажорируется
компактом \(K\). Множество всех точек пространства \(X\), инцидентных компакту
\(K\in\mathfrak S_{\kom(X)}\), мы будем обозначать через \(\mathfrak M(K)\).
При этом имеют место следующие факты:

\subsubsection\label{prop:112}
{\itshape Пусть \(K\in\mathfrak S_{\kom(X)}\). Тогда множество \(\mathfrak M(K)\)
является замкнутым} [\ref{prop:110}, \ref{pt:21}\,\eqref{eq:111}].

\subsubsection\label{prop:108}
{\itshape Пусть \(K\in\mathfrak S_{\kom(X)}\), \(x\in\mathfrak M(K)\),
\(\varepsilon\in\mathbb R^+\), и пусть \(\zeta\rightleftharpoons
\{\zeta_k\}_{k=1}^n\) "--- непустой список точек пространства \(X\), удовлетворяющий
соотношению \(\rho_{\kom(X)}(K,\zeta)<\varepsilon\). Тогда осуществим индекс
\(k\in [1,n]\), удовлетворяющий соотношению \(\rho_X(x,\zeta_k)<\varepsilon\).
}

\begin{proof}
Пусть \(\eta\) "--- непустой список точек пространства \(X\) со свойством
\begin{equation}\label{eq:102}
	\rho_{\kom(X)}(K,\eta)<\varepsilon-\rho_{\kom(X)}(K,\zeta).
\end{equation}
Тогда выполняются соотношения
\begin{flalign*}
	&& \rho_{\kom(X)}(\{x\}\cup\eta,\zeta)&\leqslant
	\rho_{\kom(X)}(\{x\}\cup\eta,\{x\}\cup K)+\rho_{\kom(X)}(K,\zeta)\\
	&& &\leqslant\rho_{\kom(X)}(K,\eta)+\rho_{\kom(X)}(K,\zeta)&
	\text{[\ref{prop:110}]}&\\
	&& &<\varepsilon,&\text{[\eqref{eq:102}]}&
\end{flalign*}
означающие осуществимость внутри \(\varepsilon\)-окрестности точки \(x\)
некоторой точки из списка \(\zeta\) [\ref{pt:22}\,\eqref{eq:001},
\ref{pt:21}\,\eqref{eq:111}].
\end{proof}

\subsubsection\label{prop:111}
{\itshape Пусть \(K\in\mathfrak S_{\kom(X)}\), \(\varepsilon\in\mathbb R^+\), и пусть
\(\zeta\rightleftharpoons\{\zeta_k\}_{k=1}^n\) "--- непустой список точек
пространства \(X\), удовлетворяющий соотношению \(\rho_{\kom(X)}(K,\zeta)<
\varepsilon\). Тогда при произвольно фиксированном \(k\in [1,n]\) осуществима точка
\(x\in\mathfrak M(K)\), удовлетворяющая соотношению \(\rho_X(x,\zeta_k)<
\varepsilon\).
}

\begin{proof}
Пусть \(\{\eta_l\}_{l=1}^{\infty}\) "--- сходящаяся в пространстве \(\kom(X)\)
к компакту \(K\) последовательность непустых списков точек пространства \(X\),
обладающая свойством
\begin{equation}\label{eq:113}
	(\forall l\geqslant 1)\qquad\rho_{\kom(X)}(\eta_l,\eta_{l+1})<2^{-l-1}.
\end{equation}
Зафиксируем номер \(m\geqslant 1\), удовлетворяющий соотношению
\begin{equation}\label{eq:116}
	2^{-m+1}<\varepsilon-\rho_{\kom(X)}(K,\zeta),
\end{equation}
и рассмотрим последовательность \(\{x_{l}\}_{l=1}^{\infty}\) точек пространства
\(X\) со следующими свойствами:
\paragraph\label{para:431} Точка \(x_1\) удовлетворяет соотношению
\(\rho_X(x_1,\zeta_k)<\rho_{\kom(X)}(K,\zeta)+2^{-m}\) и принадлежит списку \(\eta_m\).
\paragraph\label{para:432} При любом \(l\geqslant 1\) точка \(x_{l+1}\)
удовлетворяет соотношению \(\rho_X(x_{l},x_{l+1})<2^{-m-l}\) и принадлежит списку
\(\eta_{m+l}\).

\medskip\noindent
Осуществимость такой последовательности \(\{x_l\}_{l=1}^{\infty}\) гарантируется
условием \eqref{eq:113} и определением \ref{pt:21}\,\eqref{eq:111}. Для завершения
доказательства теперь остаётся заметить, что точка \(x\rightleftharpoons
\lim\limits_{l\to\infty}x_l\) удовлетворяет неравенству \(\rho_X(x,\zeta_k)<
\varepsilon\) [\ref{para:431}, \ref{para:432}, \eqref{eq:116}] и равенствам
\begin{flalign*}
	&& \rho_{\kom(X)}(\{x\}\cup K,K)&=\lim\limits_{l\to\infty}
	\rho_{\kom(X)}(\{x_{l+1}\}\cup\eta_{m+l},\eta_{m+l})\\
	&& &=0,&\text{[\ref{para:432}, \ref{pt:22}\,\eqref{eq:001},
	\ref{pt:21}\,\eqref{eq:111}]}&
\end{flalign*}
означающим выполнение соотношения \(x\in\mathfrak M(K)\).
\end{proof}

\subsubsection\label{prop:95}
{\itshape Пусть \(K_1,K_2\in\mathfrak S_{\kom(X)}\), \(x\in\mathfrak M(K_1\cup K_2)\)
и \(\varepsilon\in\mathbb R^+\). Тогда внутри \(\varepsilon\)-окрестности
точки \(x\) осуществима или точка \(x'\in\mathfrak M(K_1)\), или точка
\(x'\in\mathfrak M(K_2)\).
}

\begin{proof}
Пусть \(\zeta_1\) и \(\zeta_2\) "--- два непустых списка точек пространства \(X\),
удовлетворяющие соотношению
\begin{equation}\label{eq:1001}
	(\forall k\in\{1,2\})\qquad \rho_{\kom(X)}(K_k,\zeta_k)<\varepsilon/2.
\end{equation}
Тогда выполняется соотношение
\begin{flalign*}
	&& \rho_{\kom(X)}(K_1\cup K_2,\zeta_1\cup\zeta_2)&<\varepsilon/2,&
	\text{[\ref{prop:110}, \eqref{eq:1001}]}&
\end{flalign*}
означающее, что внутри \((\varepsilon/2)\)-окрестности точки \(x\) осуществима
или точка из списка \(\zeta_1\), или точка из списка \(\zeta_2\) [\ref{prop:108},
\ref{pt:22}\,\eqref{eq:001}]. В первом случае внутри \(\varepsilon\)-окрестности
точки \(x\) осуществима точка \(x'\in\mathfrak M(K_1)\) [\eqref{eq:1001},
\ref{prop:111}]. Во втором случае внутри \(\varepsilon\)-окрестности точки \(x\)
осуществима точка \(x'\in\mathfrak M(K_2)\) [\eqref{eq:1001}, \ref{prop:111}].
\end{proof}

\subsubsection\label{prop:115}
{\itshape Пусть \(K,K'\in\mathfrak S_{\kom(X)}\). Тогда соотношения \(K'\leqslant K\)
и \(\mathfrak M(K')\subseteq\mathfrak M(K)\) равносильны.
}

\begin{proof}
В случае \(K'\leqslant K\) факт выполнения соотношения \(\mathfrak M(K')\subseteq
\mathfrak M(K)\) тривиален [\ref{prop:109}]. В случае \(\mathfrak M(K')\subseteq
\mathfrak M(K)\) зафиксируем произвольные число \(\varepsilon\in\mathbb R^+\)
и два непустых списка \(\zeta\) и \(\zeta'\rightleftharpoons\{\zeta'_k\}_{k=1}^n\)
точек пространства \(X\) со свойствами
\begin{align}\label{eq:117}
	\rho_{\kom(X)}(K,\zeta)&<\varepsilon,\\ \label{eq:118}
	\rho_{\kom(X)}(K',\zeta')&<\varepsilon.
\end{align}
Тогда осуществим список \(\{x'_k\}_{k=1}^n\) точек множества \(\mathfrak M(K')\)
со свойством
\begin{flalign}\label{eq:119}
	&& (\forall k\in [1,n])\qquad \rho_X(x'_k,\zeta'_k)&<\varepsilon.&
	\text{[\eqref{eq:118}, \ref{prop:111}]}&
\end{flalign}
При этом для любого \(k\in [1,n]\) будет выполняться соотношение
\(x'_k\in\mathfrak M(K)\), а потому и соотношения
\begin{flalign*}
	&& \rho_{\kom(X)}(\{x'_k\}\cup\zeta,\zeta)&\leqslant
	\rho_{\kom(X)}(\{x'_k\}\cup\zeta,\{x'_k\}\cup K)+\rho_{\kom(X)}(K,\zeta)\\
	&& &\leqslant 2\rho_{\kom(X)}(K,\zeta)&\text{[\ref{prop:110}]}&\\
	&& &<2\varepsilon.&\text{[\eqref{eq:117}]}&
\end{flalign*}
Объединяя последние оценки с соотношениями \eqref{eq:119},
\ref{pt:22}\,\eqref{eq:001} и \ref{pt:21}\,\eqref{eq:111}, устанавливаем, что
при произвольно фиксированном \(k\in [1,n]\) внутри \(3\varepsilon\)-окрестности
точки \(\zeta'_k\) осуществима точка из списка \(\zeta\). В таком случае
выполняется соотношение
\begin{flalign}\label{eq:115}
	&& \rho_{\kom(X)}(\zeta'\cup\zeta,\zeta)&<3\varepsilon,&
	\text{[\ref{pt:22}\,\eqref{eq:001}, \ref{pt:21}\,\eqref{eq:111}]}&
\end{flalign}
а потому и соотношения
\begin{flalign*}
	&& \rho_{\kom(X)}(K'\cup K,K)&\leqslant\rho_{\kom(X)}(K'\cup K,\zeta'\cup\zeta)+
	\rho_{\kom(X)}(\zeta'\cup\zeta,\zeta)+\rho_{\kom(X)}(\zeta,K)\\
	&& &<\rho_{\kom(X)}(K',\zeta')+2\rho_{\kom(X)}(K,\zeta)+3\varepsilon&
	\text{[\ref{prop:110}, \eqref{eq:115}]}&\\
	&& &<6\varepsilon.&\text{[\eqref{eq:117}, \eqref{eq:118}]}&
\end{flalign*}
Ввиду произвольности произведённого ранее выбора числа \(\varepsilon\in
\mathbb R^+\), полученные оценки означают выполнение соотношения
\(K'\leqslant K\).
\end{proof}

\subsubsection\label{prop:99}
{\itshape Пусть \(K\in\mathfrak S_{\kom(X)}\), \(x\in\mathfrak S_X\)
и \(\varepsilon\in\mathbb R^+\). Тогда или множество \(\mathfrak M(K)\)
не пересекается с \(\varepsilon\)-окрестностью точки \(x\), или осуществим непустой
компакт \(K'\leqslant K\), для которого множество \(\mathfrak M(K')\) является
подмножеством \(2\varepsilon\)-окрестности точки \(x\) и надмножеством пересечения
множества \(\mathfrak M(K)\) с \(\varepsilon\)-окрестностью точки \(x\).
}

\begin{proof}
Пусть \(\{\zeta_k\}_{k=1}^{\infty}\) "--- сходящаяся в пространстве \(\kom(X)\)
к компакту \(K\) последовательность непустых списков точек пространства \(X\),
обладающая свойством
\begin{equation}\label{eq:114}
	(\forall k\geqslant 1)\qquad\rho_{\kom(X)}(\zeta_k,\zeta_{k+1})<
	2^{-k-2}\varepsilon.
\end{equation}
Заметим, что или ни одна точка из списка \(\zeta_1\) не лежит
в (\(5\varepsilon/4\))-окрестности точки \(x\), или хотя бы одна из точек этого
списка лежит в (\(3\varepsilon/2\))-окрестности указанной точки (см., например,
\cite[Гл.~2, \S~3, Теорема~20]{Ku:1973}). В первом случае множество
\(\mathfrak M(K)\) не пересекается с \(\varepsilon\)-окрестностью точки \(x\)
[\eqref{eq:114}, \ref{prop:108}].

Во втором случае рассмотрим последовательность \(\{\zeta'_k\}_{k=1}^{\infty}\)
непустых списков точек пространства \(X\) со следующими свойствами:
\paragraph\label{para:511} При любом \(k\geqslant 1\) список \(\zeta'_k\) является
подсписком списка \(\zeta_k\).
\paragraph\label{para:512} При любом \(k\geqslant 1\) любая точка списка
\(\zeta_k\), лежащая в \((1+2^{-k-1})\varepsilon\)-окрестности точки \(x\),
принадлежит списку \(\zeta'_k\), а любая точка списка \(\zeta_k\), не лежащая
в \((2-2^{-k-1})\varepsilon\)-окрестности точки \(x\), не принадлежит
списку \(\zeta'_k\).
\paragraph\label{para:513} При любом \(k\geqslant 1\) выполняется соотношение
\(\rho_{\kom(X)}(\zeta'_k,\zeta'_{k+1})<2^{-k-2}\varepsilon\).

\medskip\noindent
Осуществимость такой последовательности \(\{\zeta'_k\}_{k=1}^{\infty}\)
гарантируется вышеупомянутой теоремой \cite[Гл.~2, \S~3, Теорема~20]{Ku:1973},
соотношением \eqref{eq:114} и определением \ref{pt:21}\,\eqref{eq:111}. Проверим,
что непустой компакт \(K'\rightleftharpoons\lim\limits_{k\to\infty}\zeta'_k\)
обладает всеми требуемыми от искомого компакта свойствами.

Во-первых, выполняются равенства
\begin{flalign*}
	&& \rho_{\kom(X)}(K'\cup K, K)&=\lim\limits_{k\to\infty}
	\rho_{\kom(X)}(\zeta'_k\cup\zeta_k,\zeta_k)\\
	&& &=0,&\text{[\ref{para:511}, \ref{pt:22}\,\eqref{eq:001},
	\ref{pt:21}\,\eqref{eq:111}]}&
\end{flalign*}
а потому и соотношение \(K'\leqslant K\).

Во-вторых, внутри (\(\varepsilon/4\))-окрестности произвольно фиксированной точки
\(x'\in\mathfrak M(K')\) осуществима точка \(x''\) из списка \(\zeta'_1\)
[\ref{para:513}, \ref{prop:108}]. При этом выполняются соотношения
\begin{flalign*}
	&& \rho_X(x,x')&\leqslant\rho_X(x,x'')+\rho_X(x'',x')\\
	&& &<7\varepsilon/4+\varepsilon/4&\text{[\ref{para:511}, \ref{para:512}]}&\\
	&& &=2\varepsilon,
\end{flalign*}
означающие, что множество \(\mathfrak M(K')\) является подмножеством
\(2\varepsilon\)-окрестности точки \(x\).

В-третьих, зафиксируем произвольную точку \(x'\in\mathfrak M(K)\), лежащую внутри
\(\varepsilon\)-окрестности точки \(x\). При любом \(k\geqslant 1\) внутри
(\(2^{-k-1}\varepsilon\))-окрестности точки \(x'\) осуществима точка из списка
\(\zeta_k\) [\eqref{eq:114}, \ref{prop:108}], принадлежащая также списку
\(\zeta'_k\) [\ref{para:512}]. Тогда при любом \(k\geqslant 1\) внутри
(\(2^{-k}\varepsilon\))-окрестности точки \(x'\) осуществима точка \(x''\in
\mathfrak M(K')\) [\ref{para:513}, \ref{prop:111}], что с очевидностью
означает выполнение соотношения \(x'\in\mathfrak M(K')\) [\ref{prop:112}].
Тем самым, множество \(\mathfrak M(K')\) является надмножеством пересечения
множества \(\mathfrak M(K)\) с \(\varepsilon\)-окрестностью точки \(x\).
\end{proof}

\subsection
Имеет место следующий факт:

\subsubsection\label{prop:135}
{\itshape Пусть \(K\in \kom(\mathbf R)\). Тогда осуществимы вещественные числа
\(\sup\mathfrak M(K)\) и \(\inf\mathfrak M(K)\).
}

\begin{proof}
Пусть \(\{\zeta_k\}_{k=1}^{\infty}\) "--- сходящаяся в пространстве \(\kom(\mathbf R)\)
к компакту \(K\) последовательность непустых списков вещественных чисел, обладающая
свойством
\begin{equation}\label{eq:121}
	(\forall k\geqslant 1)\qquad
	\rho_{\kom(\mathbf R)}(\zeta_k,\zeta_{k+1})<2^{-k-1}.
\end{equation}
Рассмотрим две числовые последовательности \(\{\zeta^+_k\}_{k=1}^{\infty}\)
и \(\{\zeta^-_k\}_{k=1}^{\infty}\), соответственно, верхних и нижних точных
граней списков \(\zeta_k\). Эти последовательности удовлетворяют соотношениям
\begin{flalign*}
	&& (\forall k\geqslant 1)\qquad |\zeta^{\pm}_k-\zeta^{\pm}_{k+1}|&
	<2^{-k-1},&\text{[\eqref{eq:121}, \ref{pt:21}\,\eqref{eq:111}]}&
\end{flalign*}
а потому сходятся к некоторым вещественным числам \(\zeta^{\pm}\). Далее, при любом
\(k\geqslant 1\) внутри \(2^{-k}\)-окрестности каждой из точек \(\zeta^{\pm}_k\)
осуществима точка из множества \(\mathfrak M(K)\) [\eqref{eq:121}, \ref{prop:111}],
а потому каждая из точек \(\zeta^{\pm}\) принадлежит множеству \(\mathfrak M(K)\)
[\ref{prop:112}]. Наконец, любая точка \(x\in\mathfrak M(K)\) удовлетворяет
соотношению
\begin{flalign*}
	&& (\forall k\geqslant 1)\qquad x&\in(\zeta^-_k-2^{-k},\,\zeta^+_k+2^{-k}),&
	\text{[\eqref{eq:121}, \ref{prop:108}]}&
\end{flalign*}
из которого предельным переходом выводится соотношение \(x\in [\zeta^-,\zeta^+]\).
Тем самым, выполняются равенства
\begin{align*}
	\zeta^+&=\sup\mathfrak M(K),\\
	\zeta^-&=\inf\mathfrak M(K),
\end{align*}
означающие справедливость доказываемого утверждения.
\end{proof}

\subsection
Далее мы не будем явным образом прослеживать наиболее очевидные случаи применения
утверждений из настоящего параграфа.


\section{Вспомогательные утверждения}\label{par:3}
\subsection
В дальнейшем подмножество \(\mathfrak A\subseteq\mathfrak S_X\) носителя полного
метрического пространства \(X\) мы будем называть \emph{предкомпактным}, если
осуществим непустой компакт \(K\in\mathfrak S_{\kom(X)}\), для которого множество
\(\mathfrak M(K)\) является замыканием множества \(\mathfrak A\).

\subsection
Пусть \(f:X\to Y\) есть всюду определённое конструктивное отображение полного
метрического пространства \(X\) в полное метрическое пространство \(Y\). Пусть
также для любого предкомпактного подмножества \(\mathfrak A\subseteq\mathfrak S_X\)
его \(f\)-образ
\[
	\{y\in\mathfrak S_Y\mid (\exists x\in\mathfrak A)\quad y=f(x)\}
\]
также предкомпактен. Тогда любому непустому компакту \(K\in\mathfrak S_{\kom(X)}\)
может быть сопоставлен непустой компакт \(\Pi_f(K)\in\mathfrak S_{\kom(Y)}\),
для которого множество \(\mathfrak M(\Pi_f(K))\) является замыканием \(f\)-образа
множества \(\mathfrak M(K)\). При этом для любых двух компактов \(K,K'\in
\mathfrak S_{\kom(X)}\) из равенства \(K=K'\) с очевидностью вытекает равенство
\(\Pi_f(K)=\Pi_f(K')\) [\ref{par:2}.\ref{prop:115}]. Иначе говоря, с любым всюду
определённым конструктивным отображением \(f:X\to Y\) полных метрических
пространств, сохраняющим предкомпактность подмножеств, связано всюду определённое
конструктивное отображение \(\Pi_f:\kom(X)\to \kom(Y)\).

\subsection
Имеет место следующий тривиальный факт:

\subsubsection\label{prop:216}
{\itshape Пусть \(f:X\to Y\) "--- всюду определённое конструктивное отображение
полного метрического пространства \(X\) в полное метрическое пространство \(Y\),
\(K\in\mathfrak S_{\kom(X)}\), \(\varepsilon\in\mathbb R^+\),
и \(y\in\mathfrak S_Y\) "--- предельная точка \(f\)-образа множества
\(\mathfrak M(K)\). Тогда осуществима точка \(x\in\mathfrak M(K)\) со свойством
\(\rho_Y(f(x),y)<\varepsilon\).
}

\begin{proof}
Множество \(\mathfrak M(K)\) является сепарабельным, т.~е. осуществима плотная
в \(\mathfrak M(K)\) по метрике пространства \(X\) последовательность
\(\{x_k\}_{k=1}^{\infty}\) точек этого множества [\ref{par:2}.\ref{prop:108},
\ref{par:2}.\ref{prop:111}]. Тогда из теоремы Г.~С.~Цейтина о непрерывности
конструктивных отображений (см., например, \cite[Гл.~9, \S~2, Теорема~11]{Ku:1973},
\cite[Основная теорема]{Cei:1962}) следует, что квазиосуществим индекс
\(k\geqslant 1\), удовлетворяющий соотношению \(\rho_Y(f(x_k),y)<\varepsilon\).
Из принципа конструктивного подбора теперь немедленно вытекает, что такой
индекс \(k\geqslant 1\) является осуществимым. Последнее с очевидностью означает
справедливость доказываемого утверждения.
\end{proof}

\subsection
Пусть \(X\) "--- полное метрическое пространство, \(f:X\to\mathbf R\) "---
всюду определённая конструктивная функция, сохраняющая предкомпактность
подмножеств, и пусть \(K\in\mathfrak S_{\kom(X)}\). Тогда паре произвольно
фиксированных чисел \(n\in\mathbb Z\) и \(\varepsilon\in\mathbb R^+\)
мы сопоставляем класс \(\mathfrak K_{n,\varepsilon}(f,K)\subseteq
\mathfrak S_{\kom(X)}\) следующего вида:

\subsubsection\label{dif:1}
{\itshape Компакт \(K'\in\mathfrak S_{\kom(X)}\) считается принадлежащим классу
\(\mathfrak K_{n,\varepsilon}(f,K)\), если выполнены следующие условия:
\paragraph\label{dif:1:1} \(K'\leqslant K\).
\paragraph\label{dif:1:2} Для любых двух точек \(x,x'\in\mathfrak M(K)\)
со свойствами \(\rho_X(x,x')<2^{-n}\) и \(|f(x)-f(x')|>2\varepsilon\)
выполняется соотношение \(x\in\mathfrak M(K')\).
\paragraph\label{dif:1:3} Для любой точки \(x\in\mathfrak M(K')\) осуществима
точка \(x'\in\mathfrak M(K)\) со свойствами \(\rho_X(x,x')<2^{-n+4}\) и \(f(x')<
\sup\mathfrak M(\Pi_f(K))-\varepsilon\).
}

\bigskip
Имеет место следующий факт:

\subsubsection\label{prop:298}
{\itshape Пусть \(n\in\mathbb Z\), \(\varepsilon\in\mathbb R^+\) и \(K'\in
\mathfrak K_{n,\varepsilon}(f,K)\). Тогда осуществим компакт \(K''\in
\mathfrak K_{n+1,\varepsilon}(f,K)\) со свойством \(\rho_{\kom(X)}(K',K'')<2^{-n+5}\).
}

\begin{proof}
Пусть \(\zeta\rightleftharpoons\{\zeta_k\}_{k=1}^m\) "--- непустой список точек
пространства \(X\) со свойством
\begin{equation}\label{eq:150}
	\rho_{\kom(X)}(K,\zeta)<2^{-n-2}.
\end{equation}
Зафиксируем набор индексов \(\Gamma\subseteq [1,m]\), для которого при любом
\(k\in [1,m]\setminus\Gamma\) множество \(\mathfrak M(K')\) не пересекается
с \(2^{-n}\)-окрестностью точки \(\zeta_k\), а при любом \(k\in\Gamma\)
осуществим непустой компакт \(K'_k\leqslant K'\), для которого множество
\(\mathfrak M(K'_k)\) является подмножеством \(2^{-n+1}\)-окрестности точки
\(\zeta_k\) и надмножеством пересечения \(2^{-n}\)-окрестности точки
\(\zeta_k\) с множеством \(\mathfrak M(K')\) [\ref{par:2}.\ref{prop:99}].
При этом выполняются равенство \(K'=\bigcup_{k\in\Gamma}K'_k\)
[\ref{dif:1:1}, \ref{par:2}.\ref{prop:115}, \eqref{eq:150}] и соотношение
\begin{flalign}\label{eq:86}
	&& (\forall k\in\Gamma)\qquad \mathfrak M(\Pi_f(K'_k))&\subseteq
	\mathfrak M(\Pi_f(K)).&\text{[\ref{dif:1:1}]}&
\end{flalign}

Зафиксируем также набор индексов \(\Gamma'\subseteq\Gamma\), удовлетворяющий
соотношениям
\begin{gather}\label{eq:93}
	(\forall k\in\Gamma')\qquad\sup\mathfrak M(\Pi_f(K'_k))-
	\inf\mathfrak M(\Pi_f(K'_k))>\varepsilon,\\ \label{eq:97}
	(\forall k\in\Gamma\setminus\Gamma')\qquad\sup\mathfrak M(\Pi_f(K'_k))-
	\inf\mathfrak M(\Pi_f(K'_k))<2\varepsilon.
\end{gather}
Осуществимость такого набора \(\Gamma'\) гарантируется теоремой \cite[Гл.~2,
\S~3, Теорема~20]{Ku:1973}. Сопоставим теперь каждому значению \(k\in\Gamma\)
компакт \(K''_k\) вида
\begin{equation}\label{eq:91}
	K''_k\rightleftharpoons\left\{
		\begin{aligned}
			&K'_k&&\text{при }k\in\Gamma',\\
			&\{x_k\}&&\text{при }k\in\Gamma\setminus\Gamma',
		\end{aligned}
	\right.
\end{equation}
где \(x_k\) "--- произвольно фиксированная точка со свойствами
\begin{gather}\label{eq:85}
	x_k\in\mathfrak M(K),\\ \label{eq:94}
	f(x_k)<\sup\mathfrak M(\Pi_f(K))-\varepsilon,\\ \label{eq:92}
	(\exists x'\in\mathfrak M(K'_k))\qquad\rho_X(x_k,x')<2^{-n+4}.
\end{gather}
Осуществимость такой точки гарантируется условием \ref{dif:1:3}. Проверим,
что компакт \(K''\rightleftharpoons\bigcup_{k\in\Gamma}K''_k\) обладает всеми
требуемыми от искомого компакта свойствами.

Выполнение соотношения \(K''\leqslant K\) очевидно [\eqref{eq:91}, \eqref{eq:85}].

Пусть две точки \(x,x'\in\mathfrak M(K)\) удовлетворяют соотношениям
\(\rho_X(x,x')<2^{-n-1}\) и \(|f(x)-f(x')|>2\varepsilon\). Тогда выполняются
соотношения \(x,x'\in\mathfrak M(K')\) [\ref{dif:1:2}], а также осуществим
индекс \(k\in[1,m]\), для которого \(2^{-n}\)-окрестность точки \(\zeta_k\)
содержит каждую из точек \(x\) и \(x'\) [\eqref{eq:150}]. В таком случае
выполняются соотношения \(x,x'\in\mathfrak M(K'_k)\) и \(k\in\Gamma'\)
[\eqref{eq:97}], а потому и соотношение \(x\in\mathfrak M(K'')\) [\eqref{eq:91}].

Далее, рассмотрим произвольную точку \(x\in\mathfrak M(K'')\). Согласно утверждению
\ref{par:2}.\ref{prop:95}, осуществим индекс \(k\in\Gamma\), для которого
внутри \(2^{-n}\)-окрестности точки \(x\) осуществима точка \(x''\in
\mathfrak M(K''_k)\). При этом выполняется одно из соотношений \(k\in\Gamma'\)
или \(k\not\in\Gamma'\). В первом случае, по построению компакта \(K'_k\),
в \(2^{-n+2}\)-окрестности точки \(x''\) осуществима точка \(x'\in\mathfrak M(K)\)
со свойством
\begin{flalign*}
	&& f(x')&<\sup\mathfrak M(\Pi_f(K))-\varepsilon.&
	\text{[\eqref{eq:91}, \eqref{eq:93}, \eqref{eq:86}, \ref{prop:216}]}&
\end{flalign*}
Во втором случае осуществима точка \(x'\rightleftharpoons x''\in\mathfrak M(K)\)
с тем же свойством [\eqref{eq:91}, \eqref{eq:94}]. В обоих случаях построенная
точка \(x'\in\mathfrak M(K)\) лежит в \(2^{-n+3}\)-окрестности точки \(x\).

Наконец, по построению компактов \(K'_k\), имеет место соотношение
\begin{flalign}\label{eq:98}
	&&(\forall k\in\Gamma)\qquad \rho_{\kom(X)}(K'_k,\{\zeta_k\})&<2^{-n+2},&
	\text{[\ref{par:2}.\ref{prop:111}]}&
\end{flalign}
а потому и соотношения
\begin{flalign*}
	&& \rho_{\kom(X)}(K',K'')&=\rho_{\kom(X)}\left(\bigcup\nolimits_{k\in\Gamma}
	K'_k,\bigcup\nolimits_{k\in\Gamma} K''_k\right)\\
	&& &\leqslant\sup\limits_{k\in\Gamma}\rho_{\kom(X)}(K'_k,K''_k)&
	\text{[\ref{par:2}.\ref{prop:110}]}&\\
	&& &\leqslant\sup\limits_{k\in\Gamma\setminus\Gamma'}\left(
	\rho_{\kom(X)}(K'_k,\{\zeta_k\})+\rho_X(\zeta_k,x_k)\right)&
	\text{[\eqref{eq:91}]}&\\
	&& &<2^{-n+2}+2^{-n+1}+2^{-n+4}&\text{[\eqref{eq:98}, \eqref{eq:92}]}&\\
	&& &<2^{-n+5}.
\end{flalign*}
Тем самым, доказываемое утверждение справедливо.
\end{proof}

\subsection
Имеют место следующие два факта:

\subsubsection\label{prop:499}
{\itshape Пусть \(X\) "--- полное метрическое пространство, \(f:X\to\mathbf R\) "---
всюду определённая конструктивная функция, сохраняющая предкомпактность
подмножеств, \(K\in\mathfrak S_{\kom(X)}\) и \(\varepsilon\in\mathbb R^+\). Тогда
осуществимо число \(\delta\in\mathbb R^+\) со свойством
\begin{equation}\label{eq:800}
	(\forall x,x'\in\mathfrak M(K))\qquad \bigl(
	\bigl(f(x)>\sup\mathfrak M(\Pi_f(K))-\varepsilon/2\bigr)\land
	\bigl(\rho_X(x,x')<\delta\bigr)\bigr)\supset
	\bigl(|f(x)-f(x')|\leqslant 2\varepsilon\bigr).
\end{equation}
}

\begin{proof}
Согласно теореме \cite[Гл.~2, \S~3, Теорема~20]{Ku:1973}, выполняется хотя бы одно
из неравенств
\begin{align}\label{eq:801}
	\sup\mathfrak M(\Pi_f(K))-\inf\mathfrak M(\Pi_f(K))&<2\varepsilon,\\
	\label{eq:802}
	\sup\mathfrak M(\Pi_f(K))-\inf\mathfrak M(\Pi_f(K))&>\varepsilon.
\end{align}
В случае, когда выполняется неравенство \eqref{eq:801}, любое число
\(\delta\in\mathbb R^+\) с очевидностью обладает свойством \eqref{eq:800}.

В случае, когда выполняется неравенство \eqref{eq:802}, зафиксируем произвольные
точку \(z\in\mathfrak S_X\) и число \(m\in\mathbb Z\) со свойством
\begin{equation}\label{eq:803}
	\rho_{\kom(X)}(K,\{z\})<2^{-m+3}.
\end{equation}
При этом будет выполняться соотношение
\begin{flalign*}
	&& (\forall x,x'\in\mathfrak M(K))\qquad\rho_X(x,x')&<2^{-m+4},&
	\text{[\eqref{eq:803}, \ref{par:2}.\ref{prop:108}]}&
\end{flalign*}
а потому и соотношение \(K\in\mathfrak K_{m,\varepsilon}(f,K)\) [\ref{dif:1},
\ref{par:2}.\ref{prop:115}, \eqref{eq:802}, \ref{prop:216}]. Тогда из утверждения
\ref{prop:298} вытекает осуществимость последовательности \(\{K_k\}_{k=1}^{\infty}\)
непустых компактов со свойствами
\begin{gather}\label{eq:816}
	(\forall k\geqslant 1)\qquad K_k\in\mathfrak K_{m+k-1,\varepsilon}(f,K),\\
	\label{eq:2100}
	(\forall k\geqslant 1)\qquad \rho_{\kom(X)}(K_k,K_{k+1})<2^{-m-k+6}.
\end{gather}
Такая последовательность с очевидностью имеет в пространстве \(\kom(X)\) предел \(K'\)
со свойством
\begin{flalign}\label{eq:804}
	&& (\forall k\geqslant 1)\qquad\rho_{\kom(X)}(K_k,K')&<2^{-m-k+7}.&
	\text{[\eqref{eq:2100}]}&
\end{flalign}
Кроме того, из соотношения \eqref{eq:816} и условия \ref{dif:1:3} вытекает
осуществимость последовательности \(\{\zeta_k\}_{k=1}^{\infty}\) списков точек
множества \(\mathfrak M(K)\) со свойствами
\begin{gather}\label{eq:810}
	(\forall k\geqslant 1)\qquad \rho_{\kom(X)}(K_k,\zeta_k)<2^{-m-k+7},\\
	\label{eq:811}
	(\forall k\geqslant 1)\qquad \mathfrak M(\Pi_f(\zeta_k))\subseteq
	(-\infty,\,\sup\mathfrak M(\Pi_f(K))-\varepsilon).
\end{gather}
Применяя теперь к отображению \(\Pi_f:\kom(X)\to \kom(\mathbf R)\) теорему Г.~С.~Цейтина
о непрерывности (см., например, \cite[Гл.~9, \S~2, Теорема~11]{Ku:1973},
\cite[Основная теорема]{Cei:1962}), убеждаемся в осуществимости числа
\(n\geqslant 1\), удовлетворяющего соотношению
\begin{equation}\label{eq:808}
	(\forall K''\in\mathfrak S_{\kom(X)})\qquad
	\bigl(\rho_{\kom(X)}(K'',K')<2^{-m-n+8}\bigr)\supset
	\bigl(\rho_{\kom(\mathbf R)}(\Pi_f(K''),\Pi_f(K'))<\varepsilon/4\bigr).
\end{equation}
Тогда должно выполняться неравенство
\begin{flalign}\label{eq:62}
	&& \rho_{\kom(\mathbf R)}(\Pi_f(K_n),\Pi_f(\zeta_n))&<\varepsilon/2,&
	\text{[\eqref{eq:804}, \eqref{eq:810}, \eqref{eq:808}]}&
\end{flalign}
а потому и соотношение
\begin{flalign}\label{eq:814}
	&& \mathfrak M(\Pi_f(K_n))&\subseteq (-\infty,\,\sup\mathfrak M(\Pi_f(K))-
	\varepsilon/2).&\text{[\eqref{eq:811}, \eqref{eq:62},
	\ref{par:2}.\ref{prop:111}]}&
\end{flalign}
Тогда для любых двух точек \(x,x'\in\mathfrak M(K)\) со свойствами
\(f(x)>\sup\mathfrak M(\Pi_f(K))-\varepsilon/2\) и \(\rho_X(x,x')<2^{-m-n+1}\)
заведомо выполняется неравенство
\begin{flalign*}
	&& |f(x)-f(x')|&\leqslant 2\varepsilon.&
	\text{[\eqref{eq:814}, \eqref{eq:816}, \ref{dif:1:2}]}&
\end{flalign*}
Тем самым, число \(\delta\rightleftharpoons 2^{-m-n+1}\) обладает свойством
\eqref{eq:800}.
\end{proof}

\subsubsection\label{prop:500}
{\itshape Пусть \(X\) "--- полное метрическое пространство, \(f:X\to\mathbf R\) "---
всюду определённая конструктивная функция, сохраняющая предкомпактность подмножеств,
\(K\in\mathfrak S_{\kom(X)}\) и \(\varepsilon\in\mathbb R^+\). Тогда осуществимо число
\(\delta\in\mathbb R^+\) со свойством
\begin{equation}\label{eq:500}
	(\forall x,x'\in\mathfrak M(K))\qquad \bigl(\rho_X(x,x')<\delta\bigr)\supset
	\bigl(|f(x)-f(x')|\leqslant\varepsilon\bigr).
\end{equation}
}

\begin{proof}
Зафиксируем пространство \(X\), функцию \(f:X\to\mathbf R\) и число
\(\varepsilon\in\mathbb R^+\), после чего рассмотрим следующее вспомогательное
условие \(\mathcal Y\), налагаемое на натуральное число \(n\): "<каков бы
ни был компакт \(K\in\mathfrak S_{\kom(X)}\), удовлетворяющий соотношению
\begin{equation}\label{eq:701}
	\sup\mathfrak M(\Pi_f(K))-\inf\mathfrak M(\Pi_f(K))<
	\dfrac{n+1}{60}\,\varepsilon,
\end{equation}
осуществимо число \(\delta\in\mathbb R^+\) со свойством \eqref{eq:500}">.
Для доказательства утверждения \ref{prop:500} достаточно установить, что любое
натуральное число \(n\) удовлетворяет условию \(\mathcal Y\). Последнее будет
сделано нами на основе арифметической индукции.

Пусть число \(n\) удовлетворяет равенству \(n\eqcirc 0\). В таком случае
для произвольно фиксированного компакта \(K\in\mathfrak S_{\kom(X)}\), удовлетворяющего
соотношению \eqref{eq:701}, любое число \(\delta\in\mathbb R^+\) с очевидностью
обладает свойством \eqref{eq:500}.

Пусть число \(n\) допускает представление в виде \(n\eqcirc n'+1\), где \(n'\in
\mathbb N\) удовлетворяет условию \(\mathcal Y\). Тогда, согласно утверждению
\ref{prop:499}, для произвольно фиксированного компакта \(K\in\mathfrak S_{\kom(X)}\),
удовлетворяющего соотношению \eqref{eq:701}, осуществимо число \(\varkappa\in
\mathbb R^+\) со свойствами
\begin{gather}\label{eq:702}
	(\forall x,x'\in\mathfrak M(K))\qquad \bigl(
	\bigl(f(x)>\sup\mathfrak M(\Pi_f(K))-\varepsilon/8\bigr)\land
	\bigl(\rho_X(x,x')<\varkappa\bigr)\bigr)\supset
	\bigl(|f(x)-f(x')|\leqslant \varepsilon/2\bigr),\\ \label{eq:703}
	(\forall x,x'\in\mathfrak M(K))\qquad \bigl(
	\bigl(f(x)>\sup\mathfrak M(\Pi_f(K))-\varepsilon/60\bigr)\land
	\bigl(\rho_X(x,x')<\varkappa\bigr)\bigr)\supset
	\bigl(|f(x)-f(x')|\leqslant\varepsilon/15\bigr).
\end{gather}
Зафиксируем теперь список \(\zeta\rightleftharpoons\{\zeta_k\}_{k=1}^m\)
точек множества \(\mathfrak M(K)\) со свойством
\begin{equation}\label{eq:505}
	\rho_{\kom(X)}(K,\zeta)<\varkappa/4.
\end{equation}
Согласно утверждению \ref{par:2}.\ref{prop:99}, каждому индексу \(k\in [1,m]\)
можно сопоставить непустой компакт \(K_k\leqslant K\), для которого множество
\(\mathfrak M(K_k)\) является подмножеством \(\varkappa\)-окрестности точки
\(\zeta_k\) и надмножеством пересечения множества \(\mathfrak M(K)\)
с (\(\varkappa/2\))-окрестностью точки \(\zeta_k\). При этом с очевидностью
выполняется соотношение
\begin{equation}\label{eq:743}
	(\forall k\in [1,m])\qquad \mathfrak M(\Pi_f(K_k))\subseteq
	\mathfrak M(\Pi_f(K)).
\end{equation}

Далее, зафиксируем два набора индексов \(\Gamma\subseteq [1,m]\)
и \(\Gamma'\rightleftharpoons [1,m]\setminus\Gamma\), удовлетворяющие
соотношениям
\begin{align}\label{eq:706}
	(\forall k\in\Gamma)\qquad f(\zeta_k)&<
	\sup\mathfrak M(\Pi_f(K))-\varepsilon/12,\\ \label{eq:707}
	(\forall k\in\Gamma')\qquad f(\zeta_k)&>
	\sup\mathfrak M(\Pi_f(K))-\varepsilon/8.
\end{align}
Осуществимость таких наборов гарантируется теоремой \cite[Гл.~2, \S~3,
Теорема~20]{Ku:1973}. При этом, по построению компактов \(K_k\), выполняется
соотношение
\begin{flalign}\label{eq:709}
	&& (\forall k\in\Gamma)\,(\forall x\in\mathfrak M(K_k))\qquad
	f(x)&\leqslant\sup\mathfrak M(\Pi_f(K))-\varepsilon/60,&
	\text{[\eqref{eq:706}, \eqref{eq:703}]}&
\end{flalign}
а потому и соотношение
\begin{flalign*}
	&& (\forall k\in\Gamma)\qquad \sup\mathfrak M(\Pi_f(K_k))-
	\inf\mathfrak M(\Pi_f(K_k))&<\dfrac{n'+1}{60}\,\varepsilon.&
	\text{[\eqref{eq:709}, \eqref{eq:743}, \eqref{eq:701}]}&
\end{flalign*}
Тогда, согласно индуктивному предположению, осуществимо число
\(\delta'\in\mathbb R^+\) со свойством
\begin{equation}\label{eq:710}
	(\forall k\in\Gamma)\,(\forall x,x'\in\mathfrak M(K_k))\qquad
	\bigl(\rho_X(x,x')<\delta'\bigr)\supset\bigl(|f(x)-f(x')|\leqslant
	\varepsilon\bigr).
\end{equation}
Кроме того, по построению компактов \(K_k\), при любом \(k\in\Gamma'\)
для любых двух точек \(x,x'\in\mathfrak M(K_k)\) выполняются соотношения
\begin{flalign*}
	&& |f(x)-f(x')|&\leqslant |f(x)-f(\zeta_k)|+|f(\zeta_k)-f(x')|\\
	&& &\leqslant\varepsilon/2+\varepsilon/2&\text{[\eqref{eq:707},
	\eqref{eq:702}]}&\\
	&& &=\varepsilon.
\end{flalign*}
Тем самым, выполняется соотношение
\begin{equation}\label{eq:716}
	(\forall k\in\Gamma')\,(\forall x,x'\in\mathfrak M(K_k))\qquad
	|f(x)-f(x')|\leqslant\varepsilon.
\end{equation}

Положим теперь \(\delta\rightleftharpoons\inf(\delta',\varkappa/4)\) и заметим,
что для любых двух точек \(x,x'\in\mathfrak M(K)\) со свойством \(\rho_X(x,x')<
\delta\) осуществим индекс \(k\in [1,m]\), для которого
(\(\varkappa/2\))-окрестность точки \(\zeta_k\) содержит каждую из точек \(x\)
и \(x'\) [\eqref{eq:505}]. Тогда, по построению компактов \(K_k\), выполняются
соотношения \(x,x'\in\mathfrak M(K_k)\). При этом в каждом из случаев
\(k\in\Gamma\) и \(k\in\Gamma'\) справедлива оценка
\begin{flalign*}
	&& |f(x)-f(x')|&\leqslant\varepsilon.&
	\text{[\eqref{eq:710}, \eqref{eq:716}]}&
\end{flalign*}
Тем самым, число \(\delta\) обладает свойством \eqref{eq:500}.
\end{proof}


\section{Основной результат}
\subsection
Теперь мы сформулируем и докажем основной результат настоящей статьи.

\subsubsection\label{prop:200}
{\itshape Пусть \(X\) и \(Y\) "--- два полных метрических пространства,
а \(f:X\to Y\) "--- всюду определённое конструктивное отображение. Тогда функция
\(f\) является локально равномерно непрерывной в том и только том случае, когда
она сохраняет предкомпактность подмножеств.
}

\begin{proof}
Пусть функция \(f:X\to Y\) является локально равномерно непрерывной. Зафиксируем
произвольные компакт \(K\in\mathfrak S_{\kom(X)}\) и предкомпактное множество
\(\mathfrak A\subseteq\mathfrak S_X\) с замыканием \(\mathfrak M(K)\). Зафиксируем
также произвольную сходящуюся в пространстве \(\kom(X)\) к компакту \(K\)
последовательность \(\{\zeta_k\}_{k=1}^{\infty}\) списков точек множества
\(\mathfrak M(K)\), и условимся обозначать через \(\zeta_k'\) списки \(f\)-образов
точек из списков \(\zeta_k\). Ввиду равномерной непрерывности функции \(f\)
на множестве \(\mathfrak M(K)\), последовательность \(\{\zeta'_k\}_{k=1}^{\infty}\)
является фундаментальной относительно метрики пространства \(\kom(Y)\) и сходится
к некоторому компакту \(K'\in\mathfrak S_{\kom(Y)}\). При этом, также в силу
равномерной непрерывности функции \(f\) на множестве \(\mathfrak M(K)\),
для любых точки \(x\in\mathfrak A\) и числа \(\varepsilon\in\mathbb R^+\)
осуществима точка \(y\in\mathfrak M(K')\) со свойством
\begin{flalign*}
	&& \rho_Y(f(x),y)&<\varepsilon,&\text{[\ref{par:2}.\ref{prop:111}]}&
\end{flalign*}
а для любых точки \(y\in\mathfrak M(K')\) и числа \(\varepsilon\in\mathbb R^+\)
квазиосуществима точка \(x\in\mathfrak A\) с тем же свойством
[\ref{par:2}.\ref{prop:108}]. Таким образом, множество \(\mathfrak M(K')\)
представляет собой замыкание \(f\)-образа множества \(\mathfrak A\). Ввиду
произвольности выбора компакта \(K\in\mathfrak S_{\kom(X)}\) и множества
\(\mathfrak A\subseteq\mathfrak S_X\) с замыканием \(\mathfrak M(K)\),
сказанное означает сохранение функцией \(f\) предкомпактности подмножеств.

Пусть теперь функция \(f:X\to Y\) сохраняет предкомпактность подмножеств.
Зафиксируем произвольным образом компакт \(K\in\mathfrak S_{\kom(X)}\), число
\(\varepsilon\in\mathbb R^+\) и непустой список \(\zeta\rightleftharpoons
\{\zeta_k\}_{k=1}^n\) точек пространства \(Y\) со свойством
\begin{equation}\label{eq:201}
	\rho_{\kom(Y)}(\Pi_f(K),\zeta)<\varepsilon/3.
\end{equation}
Рассмотрим также набор \(\{\varphi_k\}_{k=1}^n\) функций \(\varphi_k:Y\to\mathbf R\)
вида
\begin{equation}\label{eq:3000}
	(\forall k\in [1,n])\,(\forall y\in\mathfrak S_Y)\qquad
	\varphi_k(y)=\rho_Y(y,\zeta_k).
\end{equation}
Каждая из функций \(\varphi_k\) с очевидностью является равномерно непрерывной,
а потому сохраняет предкомпактность подмножеств. Тогда каждая из функций
\(\varphi_k\circ f:X\to\mathbf R\) также сохраняет предкомпактность подмножеств,
а потому, согласно утверждению \ref{par:3}.\ref{prop:500}, осуществимо число
\(\delta\in\mathbb R^+\), удовлетворяющее соотношению
\begin{equation}\label{eq:202}
	(\forall k\in [1,n])\,(\forall x,x'\in\mathfrak M(K))\qquad
	(\rho_X(x,x')<\delta)\supset (|\varphi_k(f(x))-\varphi_k(f(x'))|\leqslant
	\varepsilon/3).
\end{equation}
С другой стороны, для произвольно фиксированной точки \(x\in\mathfrak M(K)\)
осуществим индекс \(k\in [1,n]\) со свойством
\begin{flalign}\label{eq:203}
	&& \rho_{Y}(f(x),\zeta_k)&<\varepsilon/3.&\text{[\eqref{eq:201}]}&
\end{flalign}
Тогда для любой точки \(x'\in\mathfrak M(K)\), удовлетворяющей условию
\(\rho_X(x,x')<\delta\), выполняются соотношения
\begin{flalign*}
	&& \rho_Y(f(x),f(x'))&\leqslant\rho_Y(f(x),\zeta_k)+\rho_Y(\zeta_k,f(x'))\\
	&& &=2\rho_Y(f(x),\zeta_k)+[\varphi_k(f(x'))-\varphi_k(f(x))]&
	\text{[\eqref{eq:3000}]}&\\
	&& &<\varepsilon.&\text{[\eqref{eq:203}, \eqref{eq:202}]}&
\end{flalign*}
Тем самым, выполняется соотношение
\[
	(\forall x,x'\in\mathfrak M(K))\qquad (\rho_X(x,x')<\delta)\supset
	(\rho_Y(f(x),f(x'))<\varepsilon).
\]
Ввиду произвольности произведённого ранее выбора компакта \(K\in\mathfrak S_{\kom(X)}\)
и числа \(\varepsilon\in\mathbb R^+\), полученный результат означает локальную
равномерную непрерывность функции \(f\).
\end{proof}

\end{document}